\documentclass[12pt]{article}
\usepackage{indentfirst, latexsym, bm, amssymb}

\begin{document}

\newtheorem{theorem}{Theorem}[section]
\newtheorem{corollary}[theorem]{Corollary}
\newtheorem{definition}[theorem]{Definition}
\newtheorem{conjecture}[theorem]{Conjecture}
\newtheorem{question}[theorem]{Question}
\newtheorem{lemma}[theorem]{Lemma}
\newtheorem{proposition}[theorem]{Proposition}
\newtheorem{example}[theorem]{Example}
\newenvironment{proof}{\noindent {\bf
Proof.}}{\rule{3mm}{3mm}\par\medskip}
\newcommand{\remark}{\medskip\par\noindent {\bf Remark.~~}}
\newcommand{\pp}{{\it p.}}
\newcommand{\de}{\em}

\newcommand{\JEC}{{\it Europ. J. Combinatorics},  }
\newcommand{\JCTB}{{\it J. Combin. Theory Ser. B.}, }
\newcommand{\JCT}{{\it J. Combin. Theory}, }
\newcommand{\JGT}{{\it J. Graph Theory}, }
\newcommand{\ComHung}{{\it Combinatorica}, }
\newcommand{\DM}{{\it Discrete Math.}, }
\newcommand{\ARS}{{\it Ars Combin.}, }
\newcommand{\SIAMDM}{{\it SIAM J. Discrete Math.}, }
\newcommand{\SIAMADM}{{\it SIAM J. Algebraic Discrete Methods}, }
\newcommand{\SIAMC}{{\it SIAM J. Comput.}, }
\newcommand{\ConAMS}{{\it Contemp. Math. AMS}, }
\newcommand{\TransAMS}{{\it Trans. Amer. Math. Soc.}, }
\newcommand{\AnDM}{{\it Ann. Discrete Math.}, }
\newcommand{\NBS}{{\it J. Res. Nat. Bur. Standards} {\rm B}, }
\newcommand{\ConNum}{{\it Congr. Numer.}, }
\newcommand{\CJM}{{\it Canad. J. Math.}, }
\newcommand{\JLMS}{{\it J. London Math. Soc.}, }
\newcommand{\PLMS}{{\it Proc. London Math. Soc.}, }
\newcommand{\PAMS}{{\it Proc. Amer. Math. Soc.}, }
\newcommand{\JCMCC}{{\it J. Combin. Math. Combin. Comput.}, }
\newcommand{\GC}{{\it Graphs Combin.}, }

\title{On the Erd\H{o}s-S\'{o}s Conjecture for Graphs on $n=k+4$ Vertices
\thanks{
This work is supported by National Natural Science
Foundation of China (No.11271256), Innovation Program of Shanghai Municipal Education Commission (No.14ZZ016) and Specialized Research Fund for the Doctoral Program of Higher Education (No.20130073110075).
\newline \indent $^{\dagger}$Correspondent author:
Xiao-Dong Zhang (Email: xiaodong@sjtu.edu.cn)}}
\author{ Long-Tu Yuan  and Xiao-Dong Zhang$^{\dagger}$   \\
{\small Department of Mathematics, and Ministry of Education }\\
{\small Key Laboratory of Scientific and Engineering Computing, }\\
{\small Shanghai Jiao Tong University} \\
{\small  800 Dongchuan Road, Shanghai, 200240, P.R. China}\\
}
\date{}
\maketitle
\begin{abstract}
  The Erd\H{o}s-S\'{o}s Conjecture states that if $G$ is a simple graph of order $n$ with average degree more than $k-2,$ then $G$ contains every tree of order $k$. In this  paper, we prove that Erd\H{o}s-S\'{o}s Conjecture is true for $n=k+4$.
\end{abstract}


\vskip 0.5cm

\section{Introduction}
The graphs considered in this paper are finite,undirected, and simple (no loops or multiple edges). Let $G=(V(G), E(G))$ be a simple graph of order $n$, where $V(G)$ is the vertex set  and $E(G)$ is the edge set with size $e(G)$. The {\it degree} of $v\in V(G)$, the number of edges incident to $v$, is denoted $d_{G}(v)$ and the set of neighbors of $v$ is denoted $N(v)$. If $u$ and $v$ in $V(G)$ are adjacent, we say that $u$ {\it hits} $v$ or $v$ {\it hits} $u$. If $u$ and $v$ are not adjacent, we say that $u$ {\it misses} $v$ or $v$ {\it misses} $u$. If $S\subseteq V(G)$, the induced subgraph of $G$ by $S$ is denoted by $G[S]$.  Denote by $D(G)$ the diameter of $G$.  In addition, $\delta(G)$, $\Delta(G)$  and $avedeg(G)=\frac{2e(H)}{|V(H)|}$  are denoted by the minimum, maximum and average degree in $V(G)$, respectively.     Let $T$ be a tree on $k$ vertices.  If there exists a injection $g:V(T)\rightarrow V(G)$ such that  $g(u)g(v)\in E(G)$  if   $uv \in E(T)$ for $u, v\in V(T)$,  we call $g$ an {\it embedding} of $T$ into $G$ and $G$ contains a copy of $T$ as a subgraph, denoted by $T\subseteq G$.  In addition, assume that $T^{\prime}\subset T$ is a proper subtree of $T$ and $g^{\prime}$ is an embedding of $T^{\prime}$ into $G$. If there exists an embedding $g: V(T)\rightarrow V(G)$ such that $g(v)=g^{\prime}(v)$ for all $v\in V(T^{\prime})$, we say that $g^{\prime}$ is $T-${\it extensible.}
%\\    Let $T$ be a tree,the diameter of $T$ is denoted $D(T)$,which is the number of the edges of the longest path of T.

In 1959, Erd\H{o}s and Gallai \cite{erdHos1959maximal} proved the following theorem.
\begin{theorem}\label{erdod1959} Let $G$ be a simple graph with $avedeg(G)>k-2$. Then $G$ contains a path of order $k$.
\end{theorem}
Based the above result,  Later Erd\H{o}s and Gallai  proposed the following  well known conjecture (for example see \cite{erdos1965})
\begin{conjecture}\label{con} Let $G$ be a simple graph with $avedeg(G)>k-2$.
Then $G$ contains every tree on $k$ vertices as a subgraph.
\end{conjecture}
  Various specific cases of Conjecture~\ref{con}  have already been proven.
   For example,
   Brandt and Dobson \cite{brandt1996erdHos} proved the conjecture for graphs having girth at least 5.  Balasubramanian and Dobson \cite{balasubramanian2007erdos} proved this conjecture for graphs without  containing $K_{2,s}$,  $s<\frac{k}{12}+1$. Li,Liu and Wang \cite{Li2000erdos} proved the conjecture for graphs whose complement has girth at least 5.    In 2003, Mclennan \cite{mclennan2005erdHos} proved the following theorem.
   \begin{theorem}\label{Mclennan-dia}
    Let $G$ be a simple graph with $avedeg(G)>k-2$. Then $G$ contains every tree of order $k$ whose diameter does not excess 4 as a subgraph.
    \end{theorem}

In 2010, Eaton and Tiner \cite{eaton2010erdos} proved the the following two theorems.
\begin{theorem}\cite{eaton2010erdos}\label{eaton-del}
 Let $G$ be a simple graph with $avedeg(G)>k-2$. If $\delta(G)\geq k-4$, then G contains every tree of order $k$  as a subgraph.
\end{theorem}
\begin{theorem}\cite{eaton2010erdos}\label{k8}
 Let $G$ be a simple graph with $avedeg(G)>k-2$. If $k\le 8$, then $G$ contains every tree of order $k$ as a subgraph.
 \end{theorem}

  In 1984, Zhou \cite{Zhou1984erdHos} proved that Conjecture~\ref{con} holds for $k=n$.   Later, Wo\'{z}niak \cite{wozniak1996erdos} proved that Conjecture~\ref{con} holds for $k=n-2$.
   \begin{theorem}\cite{wozniak1996erdos}\label{k=n-2}
 Let $G$ be a simple graph of order $b$ with $avedeg(G)>k-2$. If $k=n-2$, then
 $G$ contains every tree of order $k$  as a subgraph.
 \end{theorem}
  Recently, Tiner \cite{tiner2010erdos} proved that Conjecture~\ref{con} holds for $k=n-3$ holds.
 \begin{theorem}\cite{tiner2010erdos}\label{k=n-3}
  Let $G$ be a simple graph of order $b$ with $avedeg(G)>k-2$. If $k\ge n-3$, then
 $G$ contains every tree of order $k$  as a subgraph.
 \end{theorem}
    In this paper, we establish the following:
\begin{theorem}
\label{main}
  Let $G$ be a simple graph of order $b$ with $avedeg(G)>k-2$. If $k\ge n-4$, then
 $G$ contains every tree of order $k$  as a subgraph.
 \end{theorem}

\section{Proof of Theorem~\ref{main}}
  Let $T$ be any tree of order $k$. If $k\ge n-3$, or $k\le 8$ or the diameter of $T$ is at most 4, the assertion holds by Theorems~\ref{k=n-3},\ref{k8} and \ref{Mclennan-dia}. We only consider $k=n-4\ge 9$, $D(T)\ge 5$ and prove the assertion by the induction. Clearly the assertion holds for $n=2$. Hence  assume Theorem~\ref{main} holds for all of the graphs of order fewer than $n$  and let $G$ be a graph of order $n$.  If there exists a vertex $v$  with $d_G(v)<\lfloor\frac{k}{2}\rfloor$, then $avedeg(G-v)>k-2$ and the assertion holds  by the induction hypothesis. Further, by Theorem~\ref{eaton-del}, without loss of generality,  there exists a vertex  $z$ in $V(G)$ such that  $\lfloor\frac{k}{2}\rfloor\le d_{G}(z)=\delta(G)\leq k-5$. Moreover, assume that $e(G)=1+\lfloor\frac{1}{2}(k-2)(k+4)\rfloor$.
  Let $T$ be any tree of order $k$  with the longest path  $P=a_{0}a_{1}\ldots a_{r-1}a_{r}$ and $ N_G(a_1)\setminus \{ a_2\} =\{b_{1},\ldots,b_{s}\} $  and  $N_G(a_{r-1})\setminus \{a_{r-2}\}=\{c_{1},\ldots,c_{t}\}$.
  Since  $avedeg(G)>k-2$, we can consider the following cases: $\Delta(G)=k+3, k+2, k+1, k, k-1$.
\subsection{$\Delta(G)=k+3$ }
Let $u\in V(G)$ be such vertex that $d_{G}(u)=k+3$ and let $G^{\prime}=G-\{u,z\}$ and $T^{\prime}=T-\{a_{1},b_{1},\ldots,b_{s}\}$. Then  $e(G^{\prime})\geq e(G)-(k+3)-(k-5)+1>\frac{1}{2}(k+4)(k-2)-(k+3)-(k-5)+1=\frac{1}{2}(k^{2}-2k-2)$.
So $\emph{avedeg}(G^{\prime})>(k^{2}-2k-2)/(k+2)>k-4$ and $\mid V(T^{\prime})\mid\leq k-2$. By the induction hypothesis, $T^{\prime}\subseteq G^{\prime}$.
Let $f^{\prime}$ be an embedding of $T^{\prime}$ into $G^{\prime}$. Then let $f=f^{\prime}$ in $T^{\prime}$ and  $f(a_{1})=u$, where $X=V(G)-f^{\prime}(V(T^{\prime}))$.  Since $d_{G}(u)=k+3$, $u$ hits at least $s$ vertices in $X.$  Hence $f$ can be extended to an embedding of $T$ into $G$ or we can say that $f$ is $T-$extensible.

\subsection{$\Delta (G)=k+2$}
Let $u\in V(G)$ be such vertex that $d_{G}(u)=k+2$. Then there exists only one vertex $x\in V(G)-\{u\}$ not adjacent to $u$.  We consider two subcases: $d_{G}(x)\leq k-2$ and $d_{G}(x)\geq k-1$.
\subsubsection{$d_{G}(x)\leq k-2$}
Let $G^{\prime}=G-\{u,x\}$ and $T^{\prime}=T-\{a_{1},b_{1},\ldots,b_{s}\}$. Then $e(G^{\prime})\geq e(G)-(k+2)-(k-2)>\frac{1}{2}(k+4)(k-2)-(k+2)-(k-2)=\frac{1}{2}(k^{2}-2k-8)$. So $\emph{avedeg}(G^{\prime})>(k^{2}-2k-8)/(k+2)=k-4$ and $\mid V(T^{\prime})\mid\leq k-2$. By the induction hypothesis,  $T^{\prime}\subseteq G^{\prime}$.
Then let $f^{\prime}$ be an embedding of $T^{\prime}$ into $G^{\prime}$ and let  $f=f^{\prime}$ in $T^{\prime}$ and $f(a_{1})=u$, where $X=V(G)-f^{\prime}(V(T^{\prime}))$.  Since $d_{G}(u)=k+2$,  $u$ hits at least $s$ vertices in $X$ and $f$ is $T-$extensible.
\subsubsection{$d_{G}(x)\geq k-1$}
 We consider the following two cases.

 (A).  $x$ misses $z$. Let $G^{\prime}=G-\{u,z,x\}$ and $T^{\prime}=T-\{a_{1},b_{1},\ldots,b_{s},a_{r}\}$. Then $e(G^{\prime})\geq e(G)-(k+2)-(k-5)-(k+1)+1>\frac{1}{2}(k+4)(k-2)-(k+2)-(k-5)-(k+1)+1=\frac{1}{2}(k^{2}-4k-2)$.
   Hence $\emph{avedeg}(G^{\prime})>(k^{2}-4k-2)/(k+1)>k-5$ and $\mid V(T^{\prime})\mid\leq k-3$. By the induction hypothesis,  we have $T^{\prime}\subseteq G^{\prime}$. Since $x$ misses $z,u$ and $d_{G}(x)\geq k-1$, $x$ misses at most two vertices of $ G^{\prime}$. If $x$ hits $f^{\prime}(a_{2})$, let $f(a_{1})=x$ and $ f(a_{r})=u$. Since $d_{G}(x)\geq k-1$ and $u$ hits all vertices of $T^{\prime}$, $f$ is $T-$extensible. Hence we assume that  $x$ misses $f^{\prime}(a_{2})$. If $x$ hits $f^{\prime}(a_{r-1})$, let $f(a_{r})=x$ and $f(a_{1})=u$. Then $f$ is $T-$extensible. If $x$ misses $f^{\prime}(a_{2})$ and $f^{\prime}(a_{r-1})$, then  $x$ hits all of $V(G)-\{f^{\prime}(a_{2}),f^{\prime}(a_{r-1})\} $, because  $D(T)\geq 5$, $a_{2}$ and $a_{r-1}$ are not adjacent. Then let $f(a_{r-1})=x,f(a_{1})=u$, which implies that $f$ is $T-$extensible.

(B)  $x$ hits $z$. We consider the following two subcases.

(B.1) $d_{G}(x)>k-1$. Let $G^{\prime}=G-\{u,z,x\}$,$T^{\prime}=T-\{a_{1},b_{1},\ldots,b_{s},a_{r}\}$.
Since $x$ misses $u$ and $d_{G}(x)>k-1$,$x$ misses at most two vertices of $ G^{\prime}$, the assertion can be proven by similar to method of (A).

(B.2). $d_{G}(x)=k-1$.  Then $x$ misses 3 vertices of $V(G)\setminus\{u\}$, says $y_{1},y_{2},y_{3}$.

(a).  There exists one vertex $y_{i}$ with $1\le i\le 3$ such that $d_{G}(y_{i})=k+2$. Let $G^{\prime}=G-\{u,z,y_{i},x\}$ and $T^{\prime}=T-\{a_{1},b_{1},\ldots,b_{s},a_{r-1},c_{1},\ldots,c_{t}\}$. Then $e(G^{\prime})\geq e(G)-(k+2)-(k-5)-(k+2)-(k-1)+3+1>
\frac{1}{2}(k+4)(k-2)-(k+2)-(k-5)-(k+2)-(k-1)+3+1=\frac{1}{2}(k^{2}-6k+4)$,  because $u$ hits $z, z$ hits $x$, $u$ hits $y_{i}$, and $y_{i}$ hits $z$ by $d_{G}(y_{i})=k+2$.
Thus $\emph{avedeg}(G^{\prime})>(k^{2}-6k+4)/k>k-6$ and $\mid V(T^{\prime})\mid\leq k-4$. By the induction hypothesis,  $T^{\prime}\subseteq G^{\prime}$. Let $f(a_{1})=u$ and $f(a_{r-1})=y_{i}$. Then $f$ is $T-$extensible because $u$ and $y_{i}$ hits all the vertices of $V(T^{\prime})$, respectively.

(b).  There exists one vertex $y_{i}$ with $1\le i\le 3$ such that $d_{G}(y_{i})=k+1$.  Let $G^{\prime}=G-\{u,z,y_{i},x\}$ and $T^{\prime}=T-\{a_{1},b_{1},\ldots,b_{s},a_{r-1},c_{1},\ldots,c_{t}\}$. Then $e(G^{\prime})\geq e(G)-(k+2)-(k-5)-(k+1)-(k-1)+3>\frac{1}{2}(k+4)(k-2)-(k+2)-(k-5)-(k+1)-(k-1)+3
=\frac{1}{2}(k^{2}-6k+4)$, which implies $\emph{avedeg}(G^{\prime})>(k^{2}-6k+4)/k>k-6$ and $\mid V(T^{\prime})\mid\leq k-4$. Hence by the induction hypothesis,  $T^{\prime}\subseteq G^{\prime}$. Note that  $y_{i}$  misses at most one vertex of $G^{\prime}$. If $y_{i}$ misses $f^{\prime}(a_{2})$, let $f(a_{1})=u,f(a_{r-1})=y_{i}$; if $y_{i}$ misses $f^{\prime}(a_{r-2})$, let $f(a_{r-1})=u$,$f(a_{1})=y_{i}$. Thus $f$ is $T-$extensible.

(c).  There exists one vertex $y_{i}$ with $1\le i\le 3$ such that $d_{G}(y_{i})=k$ and $y_i$ misses $z$. Then the proof is similar to (b) and omitted.

(d). There exists one vertex $y_{i}$ with $1\le i\le 3$ such that $d_{G}(y_{i})\leq k-2$. Let $G^{\prime}=G-\{u,y_{i},x\}$ and $T^{\prime}=T-\{a_{1},b_{1},\ldots,b_{s},a_{r}\}$. Then  $e(G^{\prime})\geq e(G)-(k+2)-(k-2)-(k-1)+1>\frac{1}{2}(k+4)(k-2)-(k+2)-(k-2)-(k-1)+1=
\frac{1}{2}(k^{2}-4k-4)$, which implies $\emph{avedeg}(G^{\prime})>(k^{2}-4k-4)/(k+1)>k-5$ and $\mid V(T^{\prime})\mid\leq k-3$. Hence by the induction hypothesis, $T^{\prime}\subseteq G^{\prime}$. Similar to case (A), there exists an embedding from $T$ into $G$.

(e). $d_{G}(y_{i})=k$ and $y_i$ hits $z$ for $i\in\{1,2,3\}$;  or $d_{G}(y_{i})=k-1$ for $i\in$\{1,2,3\}.

 (e.1) $d_{T}(a_{1})+d_{T}(a_{r-1})\geq 5$.  Let $G^{\prime}=G-\{u,z,y_{1},y_{2},x\}$ and $T^{\prime}=T-\{a_{1},b_{1},\ldots,b_{s},a_{r-1},c_{1},\ldots,c_{t}\}$. Then $e(G^{\prime})\geq e(G)-(k+2)-(k-5)-(k-1)-(k-1)-(k-1)+3>\frac{1}{2}(k+4)(k-2)-(k+2)-(k-5)-(k-1)-(k-1)-(k-1)+3
 =\frac{1}{2}(k^{2}-8k+10)$ which implies $\emph{avedeg}(G^{\prime})>(k^{2}-8k+10)/(k-1)>k-7$ and $\mid V(T^{\prime})\mid\leq k-5$. Hence by the induction hypothesis, $T^{\prime}\subseteq G^{\prime}$. Moreover, $x$ misses only one vertex of $G^{\prime}$. If $x$ misses $f^{\prime}(a_{2})$, let $f(a_{1})=u,f(a_{r-1})=x$; if $x$ misses $f^{\prime}(a_{r-2})$, let $f(a_{r-1})=u,f(a_{1})=x$. In all situations, $f$ is $T-$extensible.

(e.2).  $d_{T}(a_{1})=d_{T}(a_{r-1})=2$. Let $G^{\prime}=G-\{u,z\}$ and $T^{\prime}=T-\{a_{0},a_{1}\}$. Then $e(G^{\prime})\geq e(G)-(k+2)-(k-5)+1>\frac{1}{2}(k+4)(k-2)-(k+2)-(k-5)+1=\frac{1}{2}(k^{2}-2k)$, which implies $\emph{avedeg}(G^{\prime})>(k^{2}-2k)/(k+2)>k-4$ and $\mid V(T^{\prime})\mid\leq k-2$. By the induction hypothesis,  $T^{\prime}\subseteq G^{\prime}$. Moreover, $u$  hits all vertices  of $V(G)\setminus\{x\}$ and $z$ hits $x$. Let $f(a_{1})=u $ or $z$ and $f(a_{0})=z$ or $u$. Then $f$ is $T-$extensible.

\subsection{$\Delta(G)=k+1$}
Let $u\in V(G)$ be such vertex that $d_{G}(u)=k+1$  with $u$ missing vertices $x_{1}$ and $x_{2}$. Without loss of the generality, we can assume $d_{G}(x_{1})\geq d_{G}(x_{2})$ and $d_{T}(a_{1})\geq d_{T}(a_{r-1})$.

\subsubsection{$d_{T}(a_{1})+d_{T}(a_{r-1})\geq 5$}
We consider the two cases.

(A). $x_{1}$ misses $x_{2}$.

(A.1) $d_{G}(x_{1})+d_{G}(x_{2})\leq 2k-3$.
Let $G^{\prime}=G-\{u,x_{1},x_{2}\}$ and $T^{\prime}=T-\{a_{1},b_{1},\ldots,b_{s}\}$. Then $e(G^{\prime})\geq e(G)-(k+1)-(2k-3)>\frac{1}{2}(k+4)(k-2)-(k+1)-(2k-3)=\frac{1}{2}(k^{2}-4k-4)$, which implies $\emph{avedeg}(G^{\prime})>(k^{2}-4k-4)/k+1>k-5$ and $\mid V(T^{\prime})\mid\leq k-3$. Hence by the induction hypothesis,  $T^{\prime}\subseteq G^{\prime}$. Let $f(a_{1})=u$. It is easy to see that $f$ is $T-$extensible.

(A.2). $d_{G}(x_{1})+d_{G}(x_{2})\geq 2k-2$.

(a). $d_{G}(x_{1})=k-1$ Then $d_{G}(x_{2})=k-1$ and  $x_{1}$ misses $\{u,x_{2},y_{1},y_{2}\}$. If $y_{1}, y_{2}\neq z$, let $G^{\prime}=G-\{u,z,x_{1}, x_{2}, y_{1}\}$ and $T^{\prime}=T-\{a_{1},b_{1},\ldots,b_{s}\,a_{r-1}, c_{1},\dots, c_{t}\}$. Then $e(G^{\prime})\geq e(G)-(k+1)-(k-5)-(2k-2)-(k+1)+3
%>\frac{1}{2}(k+4)(k-2)-(k+1)-(k-5)-(2k-2)-(k+1)+3=
>\frac{1}{2}(k^{2}-8k+8)$, which implies $\emph{avedeg}(G^{\prime})>(k^{2}-8k+8)/(k-1)>k-7$ and $\mid V(T^{\prime})\mid\leq k-5$. Hence  by the induction hypothesis, $T^{\prime}\subseteq G^{\prime}$.
 Note that $x_{1}$  misses only one vertex of $G^{\prime}$. If $x_{1}$ misses $f^{\prime}(a_{2})$, let $f(a_{1})=u$ and $f(a_{r-1})=x_{1}$;  if $x_{1}$ misses $f^{\prime}(a_{r-1})$, let $f(a_{r-1})=u$ and $f(a_{1})=x_{1}$. In both situations,$f$ is $T-$extensible. Now assume that $y_{1}=z$ or $y_{2}=z$. Let $G^{\prime}=G-\{u,x_{1},x_{2},y_{1},y_{2}\}$ and $T^{\prime}=T-\{a_{1},b_{1},\dots,b_{s},a_{r-1},c_{1},\dots,c_{t}\}$. Then $e(G^{\prime})\ge e(G)-(k+1)-(k-5)-(2k-2)-(k+1)+2+1
 %>\frac{1}{2}(k+4)(k-2)-(k+1)-(k-5)-(2k-2)-(k+1)+2+1=
 >\frac{1}{2}(k^{2}-8k+8)$, which implies $\emph{avedeg}(G^{\prime})>(k^{2}-8k+8)/(k-1)>k-7$ and $\mid V(T^{\prime})\mid\leq k-5$. Let $f(a_{r-1})=u$ and $f(a_{1})=x_{1}$. Then $f$ is $T-$extensible.

(b). $d_{G}(x_{1})\geq k$. Let $G^{\prime}=G-\{u,z,x_{1},x_{2}\}$ and $T^{\prime}=T-\{a_{1}, b_{1},\dots, b_{s}, a_{r-1},c_{1},\dots,\\c_{t}\}$. Then $e(G^{\prime})\geq e(G)-(k+1)-(k-5)-(2k+2)+1+2
>
%\frac{1}{2}(k+4)(k-2)-(k+1)-(k-5)-(2k+2)+1+2=
\frac{1}{2}(k^{2}-6k+2)$, which implies $\emph{avedeg}(G^{\prime})>(k^{2}-6k+2)/k>k-6$ and $\mid V(T^{\prime})\mid\leq k-4$. Hence by the induction hypothesis, $T^{\prime}\subseteq G^{\prime}$. Note that  $x_{1}$  misses at most one vertex of $G^{\prime}$. If $x_{1}$ misses $f^{\prime}(a_{2})$, let $f(a_{1})=u$ and $f(a_{r-1})=x_{1}$; if $x_{1}$ misses $f^{\prime}(a_{r-2})$, let $f(a_{r-1})=u$ and $f(a_{1})=x_{1}$. In both situations, $f$ is $T-$extensible.

(B). $x_{1}$ hits $x_{2}$.

(B.1). $d_{G}(x_{1})+d_{G}(x_{2})\leq 2k-2$. Let $G^{\prime}=G-\{u,x_{1},x_{2}\}$ and $T^{\prime}=T-\{a_{1},b_{1},\ldots,b_{s}\}$. Then $e(G^{\prime})\geq e(G)-(k+1)-(2k-2)+1>
%\frac{1}{2}(k+4)(k-2)-(k+1)-(2k-2)+1=
\frac{1}{2}(k^{2}-4k-4)$, which implies $\emph{avedeg}(G^{\prime})>(k^{2}-4k-4)/(k+1)>k-5$ and $\mid V(T^{\prime})\mid\leq k-3$. Hence by the induction hypothesis,  $T^{\prime}\subseteq G^{\prime}$. Let $f(a_{1})=u$. It is easy to see that $f$ is $T-$extensible.

(B.2). $d_{G}(x_{1})+d_{G}(x_{2})\geq 2k-1$.

(a). $d_{G}(x_{1})=k$ Then $d_{G}(x_{2})=k-1$ or $k$, and $x_{1}$ misses $u, y_{1}, y_{2}$.  If $z\neq y_{1}, y_{2}$, then let $G^{\prime}=G-\{u,z,x_{1},x_{2},y_{1}\}$ and $T^{\prime}=T-\{a_{1},b_{1},\ldots,b_{s}, a_{r-1},c_{1},\dots,c_{t}\}$. Then $e(G^{\prime})\geq e(G)-(k+1)-(k-5)-2k-(k+1)+4+1>
%\frac{1}{2}(k+4)(k-2)-(k+1)-(k-5)-2k-(k+1)+4+1=
\frac{1}{2}(k^{2}-8k+8)$, which implies $\emph{avedeg}(G^{\prime})>(k^{2}-8k+8)/(k-1)>k-7$ and $\mid V(T^{\prime})\mid\leq k-5$. Hence by the induction hypothesis,  $T^{\prime}\subseteq G^{\prime}$. Note that $x_{1}$  misses only one vertex of $G^{\prime}$. If $x_{1}$ misses $f^{\prime}(a_{2})$,let $f(a_{1})=u$ and $f(a_{r-1})=x_{1}$; if $x_{1}$ misses $f^{\prime}(a_{r-2})$, let $f(a_{r-1})=u$ and $f(a_{1})=x_{1}$. In both situations, $f$ is $T-$extensible.
If $z=y_{1}$ or $y_{2}$, then let $G^{\prime}=G-\{u,x_{1},x_{2},z\}$ and $T^{\prime}=T-\{a_{1}, b_{1}, \ldots, b_{s}, a_{r-1}, c_{1}, \dots, c_{t}\}$. Then $e(G^{\prime})\geq e(G)-(k+1)-2k-(k-5)+2>
%\frac{1}{2}(k+4)(k-2)-(k+1)-2k-(k-5)+2=
\frac{1}{2}(k^{2}-6k+4)$, which implies $\emph{avedeg}(G^{\prime})>(k^{2}-6k+4)/k>k-6$ and $\mid V(T^{\prime})\mid\leq k-4$. Hence by the induction hypothesis,  $T^{\prime}\subseteq G^{\prime}$. Note that $x_{1}$ misses only  one vertex of $G^{\prime}$.
%if $x_{1}$ misses $f^{\prime}(a_{2})$,let $f(a_{1})=u,f(a_{r-1})=x_{1}$,else let $f(a_{r-1})=u,f(a_{1})=x_{1}$.In both situations,
It is easy to see that there exists an $f$ such that $f$ is $T-$extensible.

(b). $d_{G}(x_{1})=k+1$. Let $G^{\prime}=G-\{u,x_{1},x_{2},z\}$ and $T^{\prime}=T-\{a_{1}, b_{1}, \ldots, b_{s}, a_{r-1}, c_{1},\\ \dots, c_{t}\}$. Then $e(G^{\prime})\geq e(G)-(k+1)-(k-5)-(2k+2)+2>
%\frac{1}{2}(k+4)(k-2)-(k+1)-(k-5)-(2k+2)+2=
\frac{1}{2}(k^{2}-6k)$, which implies $\emph{avedeg}(G^{\prime})>(k^{2}-6k)/k=k-6$ and $\mid V(T^{\prime})\mid\leq k-4$. Hence by the induction hypothesis, $T^{\prime}\subseteq G^{\prime}$. Note that $x_{1}$ misses at most one vertex of $G^{\prime}$. It is easy to find an embedding of $T$ into $G$.

\subsubsection{$d_{T}(a_{1})=d_{T}(a_{r-1})=2$.}
(A). There exists a vertex $v\neq u$  of degree at most $k$ such that it hits both $x_{1}$ and $x_{2}$. Let $G^{\prime}=G-\{u,v\}$ and $T^{\prime}=T-\{a_{0},a_{1}\}$. Then $e(G^{\prime})\geq e(G)-(k+1)-k+1>
%\frac{1}{2}(k+4)(k-2)-(k+1)-k+1=
\frac{1}{2}(k^{2}-2k-8)$, which implies $\emph{avedeg}(G^{\prime})>(k^{2}-2k-8)/(k+2)=k-4$ and $\mid V(T^{\prime})\mid\leq k-2$. Hence by the induction hypothesis,  $T^{\prime}\subseteq G^{\prime}$. If $f^{\prime}(a_{2})$ hits $u$, then $f(a_{1})=u$. If $f^{\prime}(a_{2})$ misses $u$, then $f^{\prime}(a_{2})$=$x_{1}$ or $x_{2}$ and let $f(a_{1})=v,f(a_{0})=u$. Thus $f$ is $T-$extensible.

(B). There exists a vertex $v\neq u$ of degree at least $k+1$ such that it hits both $x_{1}$ and $x_{2}$. Then $d_{G}(v)=k+1$ and $v$ misses $y_{1}$ and $y_{2}$. Let $G^{\prime}=G-\{u,v,z\}-\{x_{1}x_{2},y_{1}y_{2}\}$ and $T^{\prime}=T-\{a_{0},a_{1},a_{r}\}$. Then $e(G^{\prime})\geq e(G)-2(k+1)-(k-5)+1-2>
%\frac{1}{2}(k+4)(k-2)-2(k+1)-(k-5)+1-2=
\frac{1}{2}(k^{2}-4k-4)$, which implies $\emph{avedeg}(G^{\prime})>(k^{2}-4k-4)/(k+1)>k-5$ and $\mid V(T^{\prime})\mid\leq k-3$. Hence by the induction hypothesis, $T^{\prime}\subseteq G^{\prime}$. If $f^{\prime}(a_{2})=x_{1}$ or $x_{2}$, and $f^{\prime}(a_{r-1})=y_{1}$ or $y_{2}$, then let $f(a_{1})=v$ and $f(a_{r})=u$. If $f^{\prime}(a_{2})=x_{1}$ and $f^{\prime}(a_{r-1})=x_{2}$, then let $f(a_{1})=v,f^{\prime}(a_{r-1})=u,$ because $u$ hits all the neighbours of $f^{\prime}(a_{r-1})$. If $f^{\prime}(a_{2})=y_{1}$,$f^{\prime}(a_{r-1})=y_{2}$, then let $f(a_{1})=u$ and $f^{\prime}(a_{r-1})=v.$  For the rest situations, it is easy to find an embedding from $T$ into $G$.

(C). There are no vertices in $V(G)\setminus\{u\}$ hitting both $x_{1}$ and $x_{2}$, and  $x_{1}$ misses $x_{2}$.  Then  $d_{G}(x_{1})+d_{G}(x_{2})\leq k+1$. Let
$G^{\prime}=G-\{u,x_{1},x_{2}\}$ and $T^{\prime}=T-\{a_{0},a_{1}\}$. Then $e(G^{\prime})\geq e(G)-(k+1)-(k+1)>
%\frac{1}{2}(k+4)(k-2)-(k+1)-(k+1)=
\frac{1}{2}(k^{2}-2k-12)$, which implies $\emph{avedeg}(G^{\prime})>(k^{2}-2k-12)/(k+1)>k-4$  and $\mid V(T^{\prime})\mid\leq k-2$. By theorem 1.7, $T^{\prime}\subseteq G^{\prime}$. Let $f(a_{1})=u$. Then $f$ is $T-$extensible.

(D). There are no vertices in $V(G)\setminus\{u\}$ hitting both $x_{1}$ and $x_{2}$, and  $x_{1}$ hits $x_{2}$. Then  $d_{G}(x_{1})+d_{G}(x_{2})\leq k+3$. If $d_{G}(x_{1})+d_{G}(x_{2})\leq k+2$, the assertion follows from (C). Hence assume that $d_{G}(x_{1})+d_{G}(x_{2})=k+3$. Then $z$ has to hit $x_{1}$ or $x_{2}$, say that $z$ hits $x_{1}$. Let $G^{\prime}=G-\{u,z\}-\{x_{1}x_{2}\}$ and $T^{\prime}=T-\{a_{0},a_{1}\}$. Then $e(G^{\prime})\geq e(G)-(k+1)-(k-5)+1-1>
%\frac{1}{2}(k+4)(k-2)-(k+1)-(k-5)+1-1=
\frac{1}{2}(k^{2}-2k)$, which implies $\emph{avedeg}(G^{\prime})>(k^{2}-2k)/(k+2)>k-4$ and $\mid V(T^{\prime})\mid\leq k-2$. Hence by the induction hypothesis,  $T^{\prime}\subseteq G^{\prime}$. If $f^{\prime}(a_{2})$ hits $u$, let $f(a_{1})=u$; if $f^{\prime}(a_{2})=x_{1}$, let $f(a_{1})=z $ and $f(a_{0})=u$. If $f^{\prime}(a_{2})=x_{2}$ and if there is a vertex $w$ in $T^{\prime}$ such that $f^{\prime}(w)=x_{1}$,  let $f^{\prime}(w)=u$, $f(a_{1})=x_{1}$ and $f(a_{0})=z$,because $u$ hits all neighbours of $f^{\prime}(w)$; if $f^{\prime}(a_{2})=x_{2}$ and there does not   exist any vertex $w$ in $T^{\prime}$ such that $f^{\prime}(w)=x_{1}$, let $f(a_{1})=x_{1},$ and $(a_{0})=z$. In all situations, $f$ is $T-$extensible.

\subsection{$\Delta(G)=k$}
Let $u\in V(G)$ be a vertex of degree $d_{G}(u)=k$  and miss three vertices $x_{1},x_{2},x_{3}$. Denote by $S=\{x_1, x_2, x_3\}$

\subsubsection{$G[S]$ contains no edges.}
Let $G^{\prime}=G-\{u\}$ and $T^{\prime}=T-\{a_{0}\}$. Then  $e(G^{\prime})\geq e(G)-k>
%\frac{1}{2}(k+4)(k-2)-k=
\frac{1}{2}(k^{2}-8)$, which implies $\emph{avedeg}(G^{\prime})>(k^{2}-8)/(k+3)>k-3$ and $\mid V(T^{\prime})\mid\leq k-1$. By the induction hypothesis,  $T^{\prime}\subseteq G^{\prime}$. If $f^{\prime}(a_{1})$ hits $u$, let $f(a_{0})=u$; if  $f^{\prime}(a_{1})=x_{i}$, $1\le i\le 3$, let $f^{\prime}(a_{1})=u$. Since $u$ hits all neighbours of $f^{\prime}(a_{1})$,  $f$ is $T-$extensible.

\subsubsection{$G[S]$ contains exactly one edge.}\label{2.42}
Without loss of the generality,  $x_{1}$ hits $x_{2}$. We consider two cases.

(A). $d_{T}(a_{1})+d_{T}(a_{r-1})\geq5$.

(A.1).  $d_{G}(x_{1})\geq k-1$ and $d_{G}(x_{2})\geq k-1$. If $x_{3}\neq z$,let $G^{\prime}=G-\{u,z,x_{3}\}-\{x_{1}x_{2}\}$ and $T^{\prime}=T-\{a_{1},b_{1},\ldots,b_{s}\}$. Then  $e(G^{\prime})\geq e(G)-k-(k-5)-k-1>
%\frac{1}{2}(k+4)(k-2)-k-(k-5)-k-1=
\frac{1}{2}(k^{2}-4k)$, which implies $\emph{avedeg}(G^{\prime})>(k^{2}-4k)/(k+1)>k-5$ and $\mid V(T^{\prime})\mid\leq k-3$. By the induction hypothesis, $T^{\prime}\subseteq G^{\prime}$. If $f^{\prime}(a_{2})$ hits $u$, then let $f(a_{1})=u$; if $f^{\prime}(a_{2})=x_{1}$ and  $x_{2}\notin f^{\prime}(V(T^{\prime}))$, then let $f(a_{1})=x_{2}$; if $f^{\prime}(a_{2})=x_{1}$ and  $x_{2}\in f^{\prime}(V(T^{\prime}))$ and $f^{\prime}(w)=x_{2}$, then let $f^{\prime}(w)=u$, $f(a_{2})=x_{1}, $ and $f(a_{1})=x_{2}$. Hence $f$ is $T-$extensible.
On the other hand, if $x_{3}=z$, let  $G^{\prime}=G-\{u,z\}-\{x_{1}x_{2}\}$ and $T^{\prime}=T-\{a_{1},b_{1},\ldots,b_{s}\}$. Similarly, we can prove that the assertion holds.

(A.2). $d_{G}(x_{3})\geq k-1$. By (A.1), we can assume that  $d_{G}(x_{1})\leq k-2$ or $d_{G}(x_{2})\leq k-2$, say $d_{G}(x_{1})\leq k-2$.
 If $z\neq x_{1}, x_{2},$  let $G^{\prime}=G-\{u,z,x_{1},x_{2},x_{3}\}$ and $T^{\prime}=T-\{a_{1},b_{1},\ldots,b_{s},a_{r-1},c_{1},\ldots,c_{t}\}$.  Then  $e(G^{\prime})\geq e(G)-k-(k-5)-(k-2)-k-k+2+1>
 %\frac{1}{2}(k+4)(k-2)-k-(k-5)-(k-2)-k-k+2+1=
 \frac{1}{2}(k^{2}-8k+12)$, which implies $\emph{avedeg}(G^{\prime})>(k^{2}-8k+12)/(k-1)>k-7$ and $\mid V(T^{\prime})\mid\leq k-5$. Hence by the induction hypothesis,  $T^{\prime}\subseteq G^{\prime}$.
 Moreover, $x_{3}$  misses at most one vertex of $V(G^{\prime})$.  If $x_{3}$ misses $f^{\prime}(a_{2})$, let $f(a_{1})=u $ and $f(a_{r-1})=x_{3}$; if $x_{3}$ hits  $f^{\prime}(a_{2})$, let $f(a_{r-1})=u$ and $f(a_{1})=x_{3}$. then $f$ is $T-$extensible.  On the other hand, if $x_{1}=z$ or $x_{2}=z$, let  $G^{\prime}=G-\{u,x_{1},x_{2},x_{3}\}$ and $T^{\prime}=T-\{a_{1},b_{1},\ldots,b_{s},a_{r-1},\\c_{1},\ldots,c_{t}\}$. Using the same above argument, we can prove the assertion.

(A.3).  $d_{G}(x_{1})=k$ and $d_{G}(x_{2})\leq k-2$. By (A.2), we can assume that $d_{G}(x_{3})\leq k-2$. If $z\neq x_{2}, x_{3}$, let $G^{\prime}=G-\{u,z,x_{1},x_{2},x_{3}\}$ and $T^{\prime}=T-\{a_{1},b_{1},\ldots,b_{s},a_{r-1},c_{1},\ldots,c_{t}\}$. Hence  $e(G^{\prime})\geq e(G)-k-(k-5)-(k-2)-k-(k-2)+2>
%\frac{1}{2}(k+4)(k-2)-k-(k-5)-(k-2)-k-k+2=
\frac{1}{2}(k^{2}-8k+10)$, which implies $\emph{avedeg}(G^{\prime})>(k^{2}-8k+10)/(k-1)>k-7$ and $\mid V(T^{\prime})\mid\leq k-5$. By the induction hypothesis, $T^{\prime}\subseteq G^{\prime}$. Note that $x_{1}$  misses at most one vertex in  $V(G^{\prime})$. If $x_{1}$ misses $f^{\prime}(a_{2})$, let $f(a_{1})=u$ and $f(a_{r-1})=x_{1}$; if $x_{1}$ hits $f^{\prime}(a_{2})$, let $f(a_{r-1})=u$ and $f(a_{1})=x_{1}$. Hence $f$ is $T-$extensible.
On the other hand, if $x_{2}=z$ or $x_{3}=z$,let  $G^{\prime}=G-\{u,x_{1},x_{2},x_{3}\}$ and $T^{\prime}=T-\{a_{1},b_{1},\ldots,b_{s},a_{r-1},c_{1},\ldots,c_{t}\}$. By the same above argument,we can prove the assertion.

(A.4). There exists at most one vertex in $\{x_1, x_2, x_3\}$ with degree at most $k-1$.
Then there exists a vertex $u^{\prime}$ in $V(G)\setminus\{x_1, x_2, x_3,u\}$ with degree at least $k-1$. Otherwise, by $\delta(G)\le k-5$, we have $\emph{avedeg}(G)\le \frac{k+(k-1)(k-2)+(k-1)+2(k-2)+(k-5)}{k+4}\le k-2$, which is a contradiction.
Let $G^{\prime}=G-\{u,u^{\prime}\}-\{x_{1}x_{2}\}$ and $T^{\prime}=T-\{a_{1},b_{1},\ldots,b_{s}\}$. Then  $e(G^{\prime})\geq e(G)-k-k+1-1>
%\frac{1}{2}(k+4)(k-2)-k-k+1-1=
\frac{1}{2}(k^{2}-2k-8)$, which implies $\emph{avedeg}(G^{\prime})>(k^{2}-2k-8)/(k+2)=k-4$ and $\mid V(T^{\prime})\mid\leq k-2$. By the induction hypothesis, $T^{\prime}\subseteq G^{\prime}$. If $f^{\prime}(a_{2})$ hits $u$,let $f(a_{1})=u$; if $f^{\prime}(a_{2})$ misses $u$, let $f(a_{2})=u$ and $f(a_{1})=u^{\prime}$. Then $f$ is $T-$extensible.

(B). $d_{T}(a_{1})=2$ and $d_{T}(a_{r-1})=2$.  If there exists a vertex $w$ that hits both $x_{1}$ and $x_{3}$, let $G^{\prime}=G-\{u,w\}-\{x_{1}x_{2}\}$ and $T^{\prime}=T-\{a_{0},a_{1}\}$. Then  $e(G^{\prime})\geq e(G)-2k+1-1>
%\frac{1}{2}(k+4)(k-2)-2k+1-1=
\frac{1}{2}(k^{2}-2k-8)$, which implies $\emph{avedeg}(G^{\prime})>(k^{2}-2k+8)/(k+2)=k-4$ and $\mid V(T^{\prime})\mid\leq k-2$. By the induction hypothesis, $T^{\prime}\subseteq G^{\prime}$. If $f^{\prime}(a_{2})=x_{1}$ or $x_{3}$, let $f(a_{1})=w$ and $f(a_{0})=u$; if $f^{\prime}(a_{2})=x_{2}$ and $x_{1}\notin f^{\prime}(V(T^{\prime}))$, let $f(a_{1})=x_{1}$ and $f(a_{0})=w$; if $f^{\prime}(a_{2})=x_{2}$ and $x_{1}\in f^{\prime}(V(T^{\prime}))$,
 %and $f^{\prime}(v)=x_{1}$,$u$ hits all neighbours of $f^{\prime}(v)$,
 let $f^{\prime}(v)=u, f(a_{1})=x_{1}$ and $f(a_{0})=w$. In the above situations, $f$ is $T-$extensible. On the other hand, if there is no vertex hits both $x_{1}$ and $ x_{3}$, or $x_{2}$ and $x_{3}$. then  $d_{G}(x_{1})+d_{G}(x_{3})\leq k$, $d_{G}(x_{2})+d_{G}(x_{3})\leq k$.  Since $d_{G}(x_{i})\geq\lfloor\frac{k}{2}\rfloor$ and $k\geq 9$,  $d_{G}(x_{i})\leq k-2$. Hence Similar to (A.4), there exists a vertex hits $u$ with degree greater than $k-1$ and an embedding of $T$ into $G$.
%\\w hits $x_{1}$,$G^{\prime}=G-\{u,w,x_{1}x_{3}\}$,$T^{\prime}=T-\{a_{0},a_{1},a_{r-1},a_{r}\}$.So $e(G^{\prime})\geq e(G)-2k-(k+1)+2>\frac{1}{2}(k+4)(k-2)-2k-(k+1)+2=\frac{1}{2}(k^{2}-4k-2)$.$\emph{avedeg}(G^{\prime})>(k^{2}-4k-2)/k>k-6(k\geq9)$and$\mid V(T^{\prime})\mid\leq k-4$.By the induction assumption,$T^{\prime}\subseteq G^{\prime}$.
%\\w hits $x_{2}$,the same as w hits $x_{1}$
%\\w hits $x_{3}$ and misses $x_{1}$ and $x_{2}$,
\subsubsection{$G[S]$ contains exactly two edges}\label{2.43-2}
Without loss of the generality, assume that $x_{1}$ hits both $x_{2}$ and $x_{3}$. We consider the two cases.

(A).  $d_{T}(a_{1})=2$. Let $G^{\prime}=G-\{u,x_{1}\}$ and $T^{\prime}=T-\{a_{0},a_{1}\}$. Then  $e(G^{\prime})\geq e(G)-2k>
%\frac{1}{2}(k+4)(k-2)-2k=
\frac{1}{2}(k^{2}-2k-8)$, which implies $\emph{avedeg}(G^{\prime})>(k^{2}-2k-8)/(k+2)>k-4$ and $\mid V(T^{\prime})\mid\leq k-2$. By the induction hypothesis, $T^{\prime}\subseteq G^{\prime}$. If $f^{\prime}(a_{2})=x_{2}$ or $x_3$ (say $x_2$), let $f(a_{1})=x_{1}$; Moreover, if $x_{3}\notin f^{\prime}(V(T^{\prime}))$,let $f(a_{0})=x_{3}$; if $x_{3}\in f^{\prime}(V(T^{\prime}))$ and $f^{\prime}(v)=x_{3}$,
%$u$ hits all neighbours of $f^{\prime}(v)$,
let $f^{\prime}(v)=u, f(a_{1})=x_{1}, $ and $f(a_{0})=x_{3}$.
Hence,$f$ is T-extensible. If $f^{\prime}(a_{2})\neq x_{2}, x_3$, then it is easy to find an embedding from $T$ to $G$.

(B). $d_{T}(a_{1})\geq 3$.

(a). $d_{G}(x_{1})\geq k-1$. If $z\neq x_{2},x_{3}$, let $G^{\prime}=G-\{u,z,x_{1}\}$ and $T^{\prime}=T-\{a_{1},b_{1},\ldots,b_{s}\}$. Then  $e(G^{\prime})\geq e(G)-k-(k-5)-k+1>
%\frac{1}{2}(k+4)(k-2)-k-(k-5)-k+1=
\frac{1}{2}(k^{2}-4k+4)$, which implies $\emph{avedeg}(G^{\prime})>(k^{2}-4k+4)/(k+1)>k-5$ and $\mid V(T^{\prime})\mid\leq k-3$. By the induction hypothesis, $T^{\prime}\subseteq G^{\prime}$. If $f^{\prime}(a_{2})=x_{2}$ or $x_3$ (say $x_2$), let $f(a_{1})=x_{1}$. Moreover, if $x_{3}\notin f^{\prime}(V(T^{\prime}))$, let $f(a_{0})=x_{3}$; if $x_{3}\in f^{\prime}(V(T^{\prime}))$  and $f^{\prime}(v)=x_{3}$, let $f^{\prime}(v)=u, f(a_{1})=x_{1}$ and $f(a_{3})=v$, because $u$ hits all neighbours of $f^{\prime}(v)$. Hence $f$ is $T-$extensible. If $f^{\prime}(a_{2})\neq x_{2}, x_3$, it is easy to find an embedding from $T$ to $G$.  On the other hand, if $z=x_{2}$ or $x_{3}$ (say $x_2$), let $G^{\prime}=G-\{u,x_{1},x_{2}\}$, by the same  argument as (a), the assertion holds.

(b). $d_{G}(x_{1})\leq k-2$,  $d_{G}(x_{2})=k$ or $d_{G}(x_{3})=k$ (say $d_{G}(x_{2})=k$.
 Then there exists a vertex $y\in V(G)\setminus\{u, x_1, x_2, x_3\}$ such that $x_{2}$ misses $y$.  So $x_2$ misses $u, x_{3}$ and $y$ and $u$ misses $x_{3}$. By
 Case~\ref{2.42}, we can assume $y$ hits $x_{3}$. Further, by (a),  we can assume  $d_{G}(y)\leq k-2$.  If $z\neq x_{1}, y$, let $G^{\prime}=G-\{u,z,x_{2},x_{3},y\}$ and $T^{\prime}=T-\{a_{1},b_{1},\ldots,b_{s},a_{r-1},c_{1},\ldots,c_{t}\}$. Then  $e(G^{\prime})\geq e(G)-k-(k-5)-k-k-(k-2)+3>
 %\frac{1}{2}(k+4)(k-2)-k-(k-5)-k-k-(k-2)+3=
 \frac{1}{2}(k^{2}-8k+12)$, which implies $\emph{avedeg}(G^{\prime})>(k^{2}-8k+12)/(k-1)>k-7$ and $\mid V(T^{\prime})\mid\leq k-5$. By the induction hypothesis,$T^{\prime}\subseteq G^{\prime}$. Further, if $f^{\prime}(a_{2})=x_{1}$,  let $f(a_{1})=x_{2}$ and $f(a_{r-1})=u$; if $f^{\prime}(a_{r-2})=x_{1}$, let $f(a_{r-1})=x_{2}$ and $f(a_{1})=u$.
 Hence $f$ is $T-$extensible. On the other hand, if $z=y$,  let $G^{\prime}=G-\{u,x_{2},x_{3},y\}$ and $T^{\prime}=T-\{a_{1},b_{1},\ldots,b_{s},a_{r-1},c_{1},\ldots,c_{t}\}$;
 if $z=x_1$, let $G^{\prime}=G-\{u,z,x_{2},x_{3},y\}$ and $T^{\prime}=T-\{a_{1},b_{1},\ldots,b_{s},a_{r-1},c_{1},\ldots,c_{t}\}$.
 Then by the same argument, it is easy to prove that the assertion holds.

(c). $d_{G}(x_{1})\leq k-2$, $d_{G}(x_{2})=k-1$ and $d_{G}(x_{3})=k-1$.
Let $G^{\prime}=G-\{u,x_{2},x_{3}\}$ and $T^{\prime}=T-\{a_{1},b_{1},\ldots,b_{s}\}$. Then  $e(G^{\prime})\geq e(G)-k-(k-1)-(k-1)>
%\frac{1}{2}(k+4)(k-2)-k-(k-1)-(k-1)=
\frac{1}{2}(k^{2}-4k-4)$, which implies $\emph{avedeg}(G^{\prime})>(k^{2}-4k-4)/(k+1)>k-5$ and $\mid V(T^{\prime})\mid\leq k-3$. By the induction hypothesis, $T^{\prime}\subseteq G^{\prime}$. If $f^{\prime}(a_{2})=x_{1}$,  let $f(a_{1})=x_{2}$, which  $f$ is $T-$extensible. If $f^{\prime}(a_{2})\neq x_{1}$, it is easy to find an embedding from $T$ to $G$.

(d). $d_{G}(x_{1})\leq k-2$, and $d_{G}(x_{2})\le k-2$ or $d_{G}(x_{3})\le k-2$ (say $d_{G}(x_{2})\le k-2$),hence $d_{G}(x_{3})\le k-1$ by (b). Then there exists a vertex $u^{\prime}\in V(G)\setminus\{x_1, x_2, x_3, u\}$ of   degree  at least $k-1$, otherwise
$2e(G)\le (k-1)(k-2)+(k-5)+k+2(k-2)+(k-1)\le (k+4)(k-2)$ which is impossible.
  Let $G^{\prime}=G-\{u,u^{\prime},x_{1}\}$ and $T^{\prime}=T-\{a_{1},b_{1},\ldots,b_{s}\}$. Then  $e(G^{\prime})\geq e(G)-2k-(k-2)+1>
  %\frac{1}{2}(k+4)(k-2)-2k-(k-2)+1=
  \frac{1}{2}(k^{2}-4k-2)$, which implies $\emph{avedeg}(G^{\prime})>(k^{2}-4k-2)/(k+1)>k-5$ and $\mid V(T^{\prime})\mid\leq k-3$. By the induction hypothesis, $T^{\prime}\subseteq G^{\prime}$.
  Hence if $f^{\prime}(a_{2})$ hits $u$, let $f(a_{1})=u$; if $f^{\prime}(a_{2})=x_{2}$ or $x_3$ (say $x_2$), let $f^{\prime}(a_{2})=u$ and $ f(a_{1})=u^{\prime}$ since $u$ hits all the neighbours of $f^{\prime}(a_{2})$. Then $f$ is  $ T-$extensible.

\subsubsection{$G[S]$ contains exactly three edges}
(A). $d_{T}(a_{1})=2$.  If there exists an $1\le i\le 3$ (say $i=1$) such that $d_{G}(x_{1})\leq k-1$, let $G^{\prime}=G-\{u,x_{1}\}-\{x_{2}x_{3}\}$ and $T^{\prime}=T-\{a_{0},a_{1}\}$. Then  $e(G^{\prime})\geq e(G)-k-(k-1)-1>
%\frac{1}{2}(k+4)(k-2)-k-(k-1)-1=
\frac{1}{2}(k^{2}-2k-8)$, which implies $\emph{avedeg}(G^{\prime})>(k^{2}-2k-8)/(k+2)>k-4$ and $\mid V(T^{\prime})\mid\leq k-2$. By the induction hypothesis, $T^{\prime}\subseteq G^{\prime}$.  If $f^{\prime}(a_{2})=x_{2}$ or $x_3$ (say $x_2$), let $f(a_{1})=x_{1}$. Moreover, if $x_{3}\notin f^{\prime}(V(T^{\prime}))$, let $f(a_{0})=x_{3}$;  and if $x_{3}\in f^{\prime}(V(T^{\prime}))$ and $f^{\prime}(v)=x_{3}$, let $f^{\prime}(v)=u,f(a_{1})=x_{1},f(a_{0})=x_{3}$ and $f(a_{3})=v$. Hence $f$ is $T-$extensible. On the other hand, if $d_{G}(x_{1})=d_{G}(x_{2})=d_{G}(x_{3})=k$, let $G^{\prime}=G-\{u,x_{1}\}$ and $T^{\prime}=T-\{a_{0},a_{1}\}$. Then  $e(G^{\prime})\geq e(G)-2k>
%\frac{1}{2}(k+4)(k-2)-2k=
\frac{1}{2}(k^{2}-2k-8)$, which implies $\emph{avedeg}(G^{\prime})>(k^{2}-2k-8)/(k+2)=k-4$ and $\mid V(T^{\prime})\mid\leq k-2$. By the induction hypothesis,$T^{\prime}\subseteq G^{\prime}$. If $f^{\prime}(a_{2})=x_{2}$ or $x_{3}$, let $ f(a_{1})=x_{1}$; if $f^{\prime}(a_{2})\neq x_{2}, x_{3}$, let $f(a_{1})=u$. Hence $f$ is $T$-extensible.

(B). $d_{T}(a_{1})\geq3$.  If there exists an $1\le i\le 3$ (say $i=1$) such that $d_{G}(x_{1})\geq k-1$,Let $G^{\prime}=G-\{u,z,x_{1}\}-\{x_{2}x_{3}\}$. By the same argument as Case~\ref{2.43-2}.(B).(a)., the assertion holds.The rest is similar as Case~\ref{2.43-2}.(B).(d).

\subsection{$\Delta(G)=k-1$}
 Since $\Delta(G)=k-1$ and $\delta(G)\ge k-5$,  there  exist  at least four vertices of  degree  $k-1$. Otherwise $2\le 3(k-1)+k(k-2)+(k-5)=(k-2)(k+4)$, which is a contradiction.  Let $u_i$ be vertex of  $d_{G}(u_{i})=k-1$ missing  four vertices of $S_i=\{x_{i1},x_{i2},x_{i3},x_{i4}\}$  for $i=1,2,3,4$.  If there exists a vertex $u_i $ with $1\le i\le 4$  such that $G[S_i]$ contains at most one edge.
    let $G^{\prime}=G-\{u_{i}\}-E(G[S_{i}])$ and $T^{\prime}=T-\{a_{0}\}$. Then $e(G^{\prime})\geq e(G)-(k-1)-1>
   %\frac{1}{2}(k+4)(k-2)-(k-1)-1=
   \frac{1}{2}(k^{2}-8)$, which implies $\emph{avedeg}(G^{\prime})>(k^{2}-8)/(k+3)>k-3$ and $\mid V(T^{\prime})\mid\leq k-1$. By the induction hypothesis, $T^{\prime}\subseteq G^{\prime}$. If $u_{i}$ hits $f^{\prime}(a_{1})$, let $f(a_{0})=u_{i}$,  and if $u_{i}$ misses $f^{\prime}(a_{1})$, let  $f^{\prime}(a_{1})=u_{i}$. Then $f$ is $T-$extensible.
Hence we assume that  $G[S_i]$ contains  at least two edges for  $i=1,2,3,4$.

\subsubsection{$d_{T}(a_{1})\geq3, d_{T}(a_{r-1})\geq2$}
(A). $G[u_1, u_2, u_3, u_4]$ contains at least one edge, say  $u_{1}$ hits $u_{2}$.  If  $z \notin S_1=\{x_{11},x_{12},x_{13},x_{14}\},$  let $G^{\prime}=G-\{u_{1},u_{2},z\}-E(G[S_{1}])$ and $T^{\prime}=T-\{a_{1},b_{1},\ldots,b_{s}\}$. Then  $e(G^{\prime})\geq e(G)-2(k-1)-(k-5)+1-6>
%\frac{1}{2}(k+4)(k-2)-2(k-1)-(k-5)+1-6=
\frac{1}{2}(k^{2}-4k-4)$, which implies $\emph{avedeg}(G^{\prime})>(k^{2}-4k-4)/(k+1)>k-5$ and $\mid V(T^{\prime})\mid\leq k-3$. By the induction hypothesis, $T^{\prime}\subseteq G^{\prime}$. Hence if $u_{1}$ hits $f^{\prime}(a_{2})$, let $f(a_{1})=u_{1}$; and if $u_{1}$ misses $f^{\prime}(a_{2})$, let $f^{\prime}(a_{2})=u_{1} $ and $f(a_{1})=u_{2}$. Since $u_{1}$ hits all the neighbours of $f^{\prime}(a_{2})$, $f$ is  $T-$extensible.  On the other hand, if  $z \in S_1=\{x_{11},x_{12},x_{13},x_{14}\}$,  say $z=x_{11}$. Let  $G^{\prime}=G-\{u_{1},u_{2},z\}-E(G[x_{12}, x_{13}, x_{13}])$. By the same argument, the assertion holds.

(B). $G[u_{1},u_{2},u_{3},u_{4}]$ contains no edges.

(B.1). If there exist two vertices, say $u_1$ and $u_2$, in $\{u_1, u_2, u_3, u_4\}$ such that  $u_{1}$ misses $y_{1}$ and $u_{2}$ misses $y_{2}$, where $y_1\neq y_2$ and $y_1, y_2\notin \{u_1, \cdots, u_4\}$. Let $G^{\prime}=G-\{u_{1},u_{2},u_{3},u_{4}\}$ and $T^{\prime}=T-\{a_{1},b_{1},\ldots,b_{s},a_{r-1},c_{1},\ldots,c_{t}\}$. Then $e(G^{\prime})\geq e(G)-4(k-1)>
%\frac{1}{2}(k+4)(k-2)-4(k-1)=
\frac{1}{2}(k^{2}-6k)$, which implies $\emph{avedeg}(G^{\prime})>(k^{2}-6k)/k=k-6$ and $\mid V(T^{\prime})\mid\leq k-4$. By the induction hypothesis, $T^{\prime}\subseteq G^{\prime}$. Hence if $f^{\prime}(a_{2})=y_{1}$, let $f(a_{1})=u_{2}$ and $f(a_{r-1})=u_{1}$;  if $f^{\prime}(a_{2})=y_{2}$, let $f(a_{1})=u_{1}$ and $f(a_{r-1})=u_{2}$. Moreover, if $f^{\prime}(a_{r-2})=y_{1}$, let $f(a_{1})=u_{1}$ and $f(a_{r-1})=u_{2}$; and if $f^{\prime}(a_{r-2})=y_{2}$, let $f(a_{1})=u_{2}$ and $f(a_{r-1})=u_{1}$. Therefore, $f$ is  $T-$extensible.

(B.2). There exist a vertex $y\notin\{u_1, \cdots, u_4\}$ such that $y$ misses $u_1, \cdots, u_4$.  Then $G[u_1, \cdots, u_4, y]$ contains no edges.

(a). $d_{T}(a_{r-1})=2$.
  Then there exists a vertex $w$ hits  $\{u_{1},u_{2},u_{3},u_{4}\}$ and $y$. Let $G^{\prime}=G-\{u_{1},w\}$ and $T^{\prime}=T-\{a_{r-1},a_{r}\}$. Then $e(G^{\prime})\geq e(G)-2(k-1)+1)>
  %\frac{1}{2}(k+4)(k-2)-2(k-1)+1=
  \frac{1}{2}(k^{2}-2k-2)$, which implies $\emph{avedeg}(G^{\prime})>(k^{2}-2k-2)/(k+2)>k-4$ and $\mid V(T^{\prime})\mid\leq k-2$. By the induction hypothesis, $T^{\prime}\subseteq G^{\prime}$. Hence if $f^{\prime}(a_{r-2})=u_{2},u_{3},u_{4}$ or $y$, let $f(a_{r-1})=w $ and $f(a_{r})=u_{1}$; and if $f^{\prime}(a_{r-2})\neq u_{2},u_{3},u_{4}, y$, let  $f(a_{r-1})=u_{1}$ and $f(a_{r})=w$. Therefore $f$ is  $T-$extensible.

(b). $d_{T}(a_{r-1})\geq 3$.  If $z\ne y$, let $G^{\prime}=G-\{u_{1},u_{2},u_{3},u_{4},y,z\}$ and $T^{\prime}=T-\{a_{1},b_{1},\ldots,b_{s},a_{r-1},c_{1},\ldots,c_{t}\}$. Then  $e(G^{\prime})\geq e(G)-4(k-1)-(k-1)-(k-5)+4>
%\frac{1}{2}(k+4)(k-2)-4(k-1)-(k-1)-(k-5)+4=
\frac{1}{2}(k^{2}-10k+20)$, which implies $\emph{avedeg}(G^{\prime})>(k^{2}-10k+20)/(k-2)>k-8$ and $\mid V(T^{\prime})\mid\leq k-6$. By the induction hypothesis, $T^{\prime}\subseteq G^{\prime}$. Let $f(a_{1})=u_{1}$ and $f(a_{r-1})=u_{2}$. Then $f$ is $T-$extensible. On the other hand, if $z=y$, let $G^{\prime}=G-\{u_{1},u_{2},u_{3},u_{4},z\}$ and $T^{\prime}=T-\{a_{1},b_{1},\ldots,b_{s},a_{r-1},c_{1},\ldots,c_{t}\}$.
By the same argument, the assertion holds.

\subsubsection{$d_{T}(a_{1})=2,d_{T}(a_{r-1})=2.$}
(A). There exists a $1\le i\le 4$, say $i=1$,  such that $G[S_1]$ contains two or three edges. If $u_{1}$ hits  one vertex, say $u_2$, of three vertices $u_2, u_3, u_4$. Let
$G^{\prime}=G-\{u_{1},u_{2}\}-E(G[S_{1}])$ and $T^{\prime}=T-\{a_{0},a_{1}\}$. Then $e(G^{\prime})\geq e(G)-2(k-1)+1-3>
%\frac{1}{2}(k+4)(k-2)-2(k-1)+1-3=
\frac{1}{2}(k^{2}-2k-8)$, which implies $\emph{avedeg}(G^{\prime})>(k^{2}-2k-8)/(k+2)=k-4$ and $\mid V(T^{\prime})\mid\leq k-2$. By the induction hypothesis, $T^{\prime}\subseteq G^{\prime}$. Hence if $u_{1}$ hits $f^{\prime}(a_{2})$, let $f(a_{1})=u_{1}$; and if $u_{1}$ misses $f^{\prime}(a_{2})$, let $f^{\prime}(a_{2})=u_{1}$ and $f(a_{1})=u_{2}$. Since $u_{1}$ hits all the neighbours of $f^{\prime}(a_{2})$, $f$ is $T-$extensible. Therefore, we assume that $u_{1}$ misses $u_{j}$ for $j=2,3,4$.
       Then $u_{1}$ misses $x_{11}=u_{2}, x_{12}=u_{3}, x_{13}=u_{4}, x_{14}$ and $G[u_2, u_3, u_4, x_{14}]$ contains two or three edges.

(A.1). $x_{14}$ hits one vertex, say $u_{2}$,  of three vertices $u_{2},u_{3},u_{4}$. Let
$G^{\prime}=G-\{u_{1},u_{2},u_{3},u_{4}\}$ and $T^{\prime}=T-\{a_{0},a_{1},a_{r-1},a_{r}\}.$ Then  $e(G^{\prime})\geq e(G)-4(k-1)
>
%\frac{1}{2}(k+4)(k-2)-4(k-1)=
\frac{1}{2}(k^{2}-6k)$, which implies $\emph{avedeg}(G^{\prime})>(k^{2}-6k)/k=k-6$ and $\mid V(T^{\prime})\mid\leq k-4$. By the induction hypothesis, $T^{\prime}\subseteq G^{\prime}$. Since $G[u_{2},u_{3},u_{4},x_{14}]$ contains two or three edges, there exists a vertex, say $u_3$, of two vertices $u_3$,$u_4$  misses at most  one vertex,  say $y_{1}$, in $V(G)\setminus\{u_{1}, u_{2}, u_{4}, x_{14}\} $.
 Hence if $f^{\prime}(a_{2})=x_{14}$ or $y_{1}$, and $f^{\prime}(a_{r-2})=y_{1}$ or $x_{14}$, let $f(a_{1})=u_{2}$ or $u_{1}$ and $f(a_{r-1})=u_{1}$ or $u_{2}$, then  $f$ is $T-$extensible. For the rest cases, it is easy to find an embedding from $T$ to $G$.

(A.2). $x_{14}$ misses three vertices $u_{2},u_{3},u_{4}$. Then $G[u_{2},u_{3},u_{4}]$ contains two or three edges. We can assume that  $u_{2}$ hits $u_{3}$ and $u_{4}$.
If $u_{3}$ misses $u_{4}$, $u_{3}$ misses at most one vertex, says $y_{1}$, in $V(G)\setminus\{u_1, u_2, u_4, x_{14}\}$. Then let $G^{\prime}=G-\{u_{1},x_{14},u_{3},u_{4}\}$ and $T^{\prime}=T-\{a_{0},a_{1},a_{r-1},a_{r}\}$. By the similar argument as Case (A.1), the assertion holds. Hence we can assume that $u_{3}$ hits $u_{4}$ and $u_{3}$ misses $z_{1},z_{2},u_{1},x_{14}$. Let $G^{\prime}=G-\{u_{1}, x_{14}, u_{3}, u_{4}\}-\{z_{1}z_{2}\}$ and $T^{\prime}=T-\{a_{0},a_{1},a_{r-1},a_{r}\}$. Then  $e(G^{\prime})\geq e(G)-4(k-1)+1-1>
%\frac{1}{2}(k+4)(k-2)-4(k-1)+1-1=
\frac{1}{2}(k^{2}-6k)$, which implies $\emph{avedeg}(G^{\prime})>(k^{2}-6k)/k=k-6$ and $\mid V(T^{\prime})\mid\leq k-4$. By the induction hypothesis, $T^{\prime}\subseteq G^{\prime}$. Hence if $f^{\prime}(a_{2})=z_{1}$ or $z_{2}$, and $f^{\prime}(a_{r-2})=z_{2}$ or $z_{1}$, let $f^{\prime}(a_{2})=u_{3}$, $f(a_{1})=u_{4}$, $f(a_{r-1})=u_{1}$. Therefore $f$ is $T-$ extensible. If  $f^{\prime}(a_{2})=z_{1}$ or $z_{2}$, and $f^{\prime}(a_{r-2})=u_{2}$,let $f(a_{1})=u_{1}$, $f(a_{r-1})=u_{4}$. Therefore $f$ is $T-$ extensible. For the rest cases, it is easy to find an embedding from $T$ to $G$.

(B).  There exists a $1\le i\le 4$, say $i=1$,  such that $G[S_1]$ contains exactly four edges.

(B.1). There exists a vertex, say $x_{11}$,  of degree 3 in $G[S_1]$ and $\mid E(G[S_1])\mid\leq 5$.  Then $x_{11}$ hits $x_{12},x_{13}$ and $x_{14}$. Let $G^{\prime}=G-\{u_{1},x_{11}\}-\{E(G[x_{12},x_{13},x_{14}])\}$ and $T^{\prime}=T-\{a_{0},a_{1}\}$. Then  $e(G^{\prime})\geq e(G)-2(k-1)-2>
%\frac{1}{2}(k+4)(k-2)-2(k-1)-1=
\frac{1}{2}(k^{2}-2k-8)$, which implies $\emph{avedeg}(G^{\prime})>(k^{2}-2k-8)/(k+2)>k-4$ and $\mid V(T^{\prime})\mid\leq k-2$. By the induction hypothesis, $T^{\prime}\subseteq G^{\prime}$.
Hence if $u_{1}$ hits $f^{\prime}(a_{2})$, let $f(a_{1})=u_{1}$, which implies that $f$ is $T-$extensible.  If $u_{1}$ misses $f^{\prime}(a_{2})$ and $f^{\prime}(a_{2})=x_{12}, $ let $f(a_{1})=x_{11}$.  Moreover, if $x_{13}$ or $x_{14} \notin f^{\prime}(V(T^{\prime}))$, then  let $f(a_{0})=x_{13}$ or $x_{14}$. Then $f$ is $T-$extensible.  If $x_{13}$ and $x_{14} \in f^{\prime}(V(T^{\prime})),$ $f^{\prime}(w)=x_{13}$ or $x_{14}$,
let $f^{\prime}(w)=u_{1},f(a_{0})=x_{13}$ or $x_{14}$. Then $f$ is $T-$extensible. For the rest cases, it is easy to find an embedding from $T$ to $G$.

(B.2). The degree of every vertex in $G[S_1]$ is two. We assume that $x_{11}$ hits $x_{12}, x_{12}$ hits $x_{13}, x_{13}$ hits $x_{14}$, $x_{14}$ hits $x_{11}.$

(a). $u_{1}$ hits all vertices of $\{u_2,u_{3},u_{4}\}$.

(a.1).
%\\(\uppercase\expandafter{\romannumeral1})
There exists a vertex $u_{i},$ say $u_2$,  in $\{u_2, u_3, u_4\}$ which  misses $x_{11},x_{12},x_{13}$ and $x_{14}.$  Let $G^{\prime}=G-\{u_{1},u_{2},x_{11},x_{12}\}-\{x_{13}x_{14}\}$ and $T^{\prime}=T-\{a_{0},a_{1},a_{r-1},a_{r}\}$. Then  $e(G^{\prime})\geq e(G)-4(k-1)+1>
%\frac{1}{2}(k+4)(k-2)-4(k-1)+1=
\frac{1}{2}(k^{2}-6k+2)$, which implies $\emph{avedeg}(G^{\prime})>(k^{2}-6k+2)/k>k-6$ and $\mid V(T^{\prime})\mid\leq k-4$. By the induction hypothesis, $T^{\prime}\subseteq G^{\prime}$.
If $f^{\prime}(a_{2})=x_{13}, f^{\prime}(a_{r-2})=x_{14}$, let  $f(a_{1})=x_{12}, f(a_{0})=x_{11}, f(a_{r-2})=u_{1}, f(a_{r-1})=u_{2}$. Hence $f$ is $T-$extensible. For the rest cases, similarly, it is easy to  find an embedding from $T$ to $G$.

%\\(\uppercase\expandafter{\romannumeral3})

(a.2). There exists a vertex, say $u_{2}$, in $\{u_2, u_3, u_4\}$ such that it  hits at least two vertices of $\{x_{11},x_{12},x_{13},x_{14}\}$, say $u_2$ hits $x_{11}$ and $x_{13}$, or $u_2$ hits $x_{11}$ and $x_{12}$.

%\\(\romannumeral1)
 If  $u_{2}$ hits $x_{11}$ and $x_{13}$, let $G^{\prime}=G-\{u_{1},u_2\}-\{x_{11}x_{12},x_{12}x_{13},x_{13}x_{14}\}$ and $T^{\prime}=T-\{a_{0},a_{1},a_{r-1},a_{r}\}$. Then  $e(G^{\prime})\geq e(G)-2(k-1)+1-3>
 %\frac{1}{2}(k+4)(k-2)-2(k-1)+1-3=
 \frac{1}{2}(k^{2}-2k-8)$, which implies $\emph{avedeg}(G^{\prime})>k-4$ and $\mid V(T^{\prime})\mid\leq k-2$. By the induction hypothesis, $T^{\prime}\subseteq G^{\prime}$. Hence if $ f^{\prime}(a_{2})=x_{11}$ or $x_{13}$, let $f(a_{1})=u_{2}$; if $f^{\prime}(a_{2})=x_{12}$, let $f(a_{2})=u_{1}$ and $f(a_{1})=u_{3}$; if $f^{\prime}(a_{2})=x_{14}$ and  $x_{13}\notin f^{\prime}(V(T^{\prime}))$,let $f(a_{1})=x_{13}$ and $f(a_{0})=u_{2}$;  if $f^{\prime}(a_{2})=x_{14}$ and $x_{13}\in f^{\prime}(V(T^{\prime}))$, let  $f(v)=u_{1},f(a_{1})=x_{13},f(a_{0})=u_{2}$, because there is a vertex $v$, $f^{\prime}(v)=x_{13}$ and $u_{1}$ hits all the neighbours of $f^{\prime}(v)$. Therefore $f$ is $T-$extensible.

  If  $u_{2}$ hits $x_{11}$ and $x_{12}$, let $G^{\prime}=G-\{u_{1},u_{2}\}-\{x_{12}x_{13},x_{13}x_{14},x_{11}x_{14}\}$ and $T^{\prime}=T-\{a_{0},a_{1}\}$. Then  $e(G^{\prime})\geq e(G)-2(k-1)+1-3>
  %\frac{1}{2}(k+4)(k-2)-2(k-1)+1-3=
  \frac{1}{2}(k^{2}-2k-8)$, which implies $\emph{avedeg}(G^{\prime})>k-4$ and $\mid V(T^{\prime})\mid\leq k-2$. By the induction hypothesis, $T^{\prime}\subseteq G^{\prime}$.  Hence if $f^{\prime}(a_{2})=x_{11}$ or $x_{12}$, let $f(a_{1})=u_{2}$; if $f^{\prime}(a_{2})=x_{13}$ or $x_{14}$, let $f(a_{2})=u_{1},f(a_{1})=u_{2}$, because $u_{1}$ hits all the neighbours of $f^{\prime}(a_{2})$. Therefore $f$ is $T-$extensible.

%(\uppercase\expandafter{\romannumeral2})
(a.3). $u_{i}$ hits exactly one vertex of $\{x_{11},x_{12},x_{13},x_{14}\}$ for $i=2, 3, 4.$

%\\(\uppercase\expandafter{\romannumeral2}.1)
(i).  There exist two vertices of $\{u_{2},u_{3},u_{4}\}$ such that they hit the same vertex in $\{x_{11},x_{12}, x_{13}, x_{14}\}$, says both $u_{2} $ and $u_{3}$ hit $x_{14}$.

   If  $u_{2}$ and $u_3$ misses the same vertices, say, $\{x_{11}, x_{12}, x_{13}, y\},u_{3}$, then $u_2$ hits $u_3$. Further, if $G[x_{11}, x_{12}, x_{13}, y]$ contains at most three edges or has a vertex of degree 3, the assertion follows from Case 2.5.2.(A) or Case 2.5.2.(B.1). Therefore we can assume that  $y$ hits  both $x_{11}$ and  $x_{13}$. Let $G^{\prime}=G-\{u_{2},u_{3},x_{11},x_{12}\}-\{x_{13}y\}$ and
   $T^{\prime}=T-\{a_{0},a_{1},a_{r-1},a_{r}\}$. The assertion follows from Case 2.5.2. (B.2).(a.1).

If $u_{2}$ misses $\{x_{11}, x_{12}, x_{13},y_1\}$ and  $u_{3}$ misses $\{x_{11}, x_{12}, x_{13},y_2\}$
 with $y_1\neq y_2$, let
 $G^{\prime}=G-\{u_{1},u_{2},u_{3},x_{14}\}-\{x_{11}x_{12},x_{12}x_{13}\}$ and $T^{\prime}=T-\{a_{0},a_{1},a_{r-1},a_{r}\}$. Then  $e(G^{\prime})\geq e(G)-4(k-1)+4-2>
 %\frac{1}{2}(k+4)(k-2)-4(k-1)+4-2=
 \frac{1}{2}(k^{2}-6k+4)$, which implies $\emph{avedeg}(G^{\prime})>k-6$ and $\mid V(T^{\prime})\mid\leq k-4$. By the induction hypothesis, $T^{\prime}\subseteq G^{\prime}$.   Hence if $f^{\prime}(a_{2})=x_{11}$ or $x_{13}$, let $f(a_{1})=x_{14}$,$f(a_{0})=u_{3}$ or $u_{2}$.if $f^{\prime}(a_{2})=x_{12}$, let $f(a_{2})=u_{1}$,$f(a_{1})=u_{3}$ or $u_{2}$.if $f^{\prime}(a_{2})=y_{1}$ or $y_{2}$, let $f(a_{1})=u_{3}$ or $u_{2}$. which implies $f$ is $T-$extensible. For the rest cases, it is easy to find an embedding from $T$ to $G$.
 %$f^{\prime}(a_{r-2})=u_{1},f(a_{r-1})=u_{2}$.2.If $f^{\prime}(a_{2})=x_{1},f^{\prime}(a_{r-2})=x_{3}$,then let $f(a_{r-1})=x_{4},f(a_{r})=u_{3},f^{\prime}(a_{r-2})=u_{1},f(a_{r-1})=u_{2}$.3.If $f^{\prime}(a_{2})=x_{1},x_{2}$ or $x_{3}$,$f^{\prime}(a_{r-2})=y_{1}$ or $y_{2}$,then let $f^{\prime}(a_{2})=u_{1},f(a_{1})=u_{3}$ or $u_{2},f(a_{r-1})=u_{3}$ or $u_{2},f(a_{r-1})=u_{2}$.4.$f^{\prime}(a_{2})=y_{1}$,$f^{\prime}(a_{r-2})=y_{2}$,then let $f(a_{1})=u_{3}$,$f(a_{r-1})=u_{2}.$Other situations are much easier.

%(\uppercase\expandafter{\romannumeral2}.2)

 (ii). $\{u_{2},u_{3},u_{4}\}$ hits the different vertices of $\{x_{11},x_{12},x_{13},x_{14}\}$.  Without loss of generality, we assume that $u_{2}$ hits $x_{11}$ and $u_{3}$ hits $x_{13}$,$u_{2}$ misses $y_{1}$ and $u_{3}$ misses $y_{2}$.
 %$u_{2}$ misses $\{x_{2},x_{3},x_{4},y_{1}\},u_{3}$ misses $\{x_{1},x_{2},x_{4},y_{2}\}$
Let $G^{\prime}=G-\{u_{1},u_{2},u_{3},x_{13}\}-\{x_{11}x_{12},x_{11}x_{14}\}$ and $T^{\prime}=T-\{a_{0},a_{1},a_{r-1},a_{r}\}$. Then  $e(G^{\prime})\geq e(G)-4(k-1)+3+0-2>
%\frac{1}{2}(k+4)(k-2)-4(k-1)+3+0-2=
\frac{1}{2}(k^{2}-6k+2)$, which implies $\emph{avedeg}(G^{\prime})>k-6$ and $\mid V(T^{\prime})\mid\leq k-4$. By the induction hypothesis, $T^{\prime}\subseteq G^{\prime}$.
%If $y_{1}=y_{2}=y:$1.
Hence if $f^{\prime}(a_{2})=x_{12}$ or $x_{14}$, let $f(a_{1})=x_{13}$ and $f(a_{0})=u_{3}$,or $f(a_{2})=u_{1}$ and $f(a_{1})=u_{2}$, if $f^{\prime}(a_{2})=y_{1}$ or $y_{2}$, let $f(a_{1})=u_{1}$, if $f^{\prime}(a_{2})=x_{11}$, let $f(a_{1})=u_{2}$, Therefore $f$ is $T-$extensible. For the rest cases, by the same argument, it is easy to find an embedding from $T$ to $G$.

%$f^{\prime}(a_{2})=y$,$f^{\prime}(a_{r-2})=x_{1}$,then let $f(a_{1})=u_{1}$,$f(a_{r-1})=u_{2}$.3.If $f^{\prime}(a_{2})=y$,$f^{\prime}(a_{r-2})=x_{2}$ or $x_{4}$,then let $f(a_{1})=u_{1}$,$f(a_{r-1})=x_{3},f(a_{r})=u_{3}$.Other situation is much easier.If $y_{1}\neq y_{2}$,add the case $f^{\prime}(a_{2})=y_{1}$,$f^{\prime}(a_{r-2})=y_{2}$,then let $f(a_{1})=u_{3},f(a_{r-1})=u_{2}$.All of the above situations are T-extensible.

(b). $u_{1}$ hits one or two vertices of $\{u_{2}, u_{3}, u_{4}\}$.
 Without loss of the generality, we assume that $u_{1}$ hits $u_{2}$  and $u_{1}$ misses $u_{4}$.  Then $u_4\in \{x_{11}, x_{12}, x_{13}, x_{14}\}$, say $u_4=x_{14}$,$u_4$ misses $u_1$,$x_{12}$,$z_1$,$z_2$.

%(b.1). $u_{2}$ hits $u_{4}$.
 If$ u_{2}\neq z_{1},z_{2}$, let $G^{\prime}=G-\{u_{1},u_{2},u_{4},x_{12}\}-\{z_{1}z_{2}\}$ and $T^{\prime}=T-\{a_{0},a_{1},a_{r-1},a_{r}\}$. Then  $e(G^{\prime})\geq e(G)-4(k-1)+1-1>
%\frac{1}{2}(k+4)(k-2)-4(k-1)+2-1=
\frac{1}{2}(k^{2}-6k)$, which implies $\emph{avedeg}(G^{\prime})
%>(k^{2}-6k+2)/k
>k-6$ and $\mid V(T^{\prime})\mid\leq k-4$. By the induction hypothesis, $T^{\prime}\subseteq G^{\prime}$. Hence if $f^{\prime}(a_{2})=x_{11}$ and $f^{\prime}(a_{r-2})=x_{13}$, let $f(a_{1})=u_{4},f(a_{r-2})=u_{1}$ and $f(a_{r-1})=u_{2}$. Therefore $f$ is T-extensible. For the rest cases, it is easy to find an embedding from $T$ to $G$. If $u_{2}=z_{1}$ or $z_{2}$, say $u_{2}=z_{1}$, let $G^{\prime}=G-\{u_{1},u_{2},u_{4},x_{12}\}$  and $T^{\prime}=T-\{a_{0},a_{1},a_{r-1},a_{r}\}$.This situation is much easier than the above case.

%(b.2). $u_{2}$ misses $u_{4}$. ,then let $y_{1}=u_{2}$.
%\\$G^{\prime}=G-\{u_{1},u_{2},u_{4},x_{2}\}$,$T^{\prime}=T-\{a_{0},a_{1},a_{r-1},a_{r}\}$.So $e(G^{\prime})\geq e(G)-4(k-1)+1>\frac{1}{2}(k+4)(k-2)-4(k-1)+1=\frac{1}{2}(k^{2}-6k+2)$.$\emph{avedeg}(G^{\prime})>(k^{2}-6k+2)/k>k-6$ and $\mid V(T^{\prime})\mid\leq k-4$.By the induction assumption,$T^{\prime}\subseteq G^{\prime}$.If $f^{\prime}(a_{2})=x_{1}$,$f^{\prime}(a_{r-2})=x_{3}$,then let $f(a_{1})=u_{4},f^{\prime}(a_{r-2})=u_{1},f(a_{r-1})=u_{2}$.If $f^{\prime}(a_{2})=x_{1}$ or $x_{3},f^{\prime}(a_{r-2})=y_{2}$,then let $f(a_{1})=u_{4},f(a_{r-1})=u_{1}$.$f$ is T-extensible,and other situations are much easier.

(c).  $u_{1}$ misses all vertices of $\{u_2, u_3, u_4\}$.
 Without loss of generality, we assume  $u_{2}=x_{11},$ $u_3=x_{12},$ $u_4=x_{13}$.
  Let $u_{2}$ miss  $\{u_{1}, x_{13}, y_1, y_2\}$.
  If $G[u_{1}, x_{13}, y_1, y_2]$ contains two, or three edges, or a vertex of degree 3, the assertion follows from Case 2.5.2 (A). and Case 2.5.2 (B.1).
  Hence we assume that $u_1$ hits $y_1$, $y_1$ hits $u_4=x_{13}$, $u_4$ hits $y_2$ and $y_2$ hits $u_1$.  Hence the assertion follows from Case 2.5.2. (B.2). (a)  and Case 2.5.2. (B.2).(b).

(C). There exists a $1\le i\le 4$, say $i=1$, such that  $G[x_{11},x_{12},x_{13},x_{14}]$ contains five edges. Then we assume that
 $x_{11}$ hits $x_{12},x_{13}$ and $x_{14}. $
  Let $G^{\prime}=G-\{u_{1},x_{11}\}-\{E(G[x_{12},x_{13},x_{14}])\}$ and $T^{\prime}=T-\{a_{0},a_{1}\}$. The assertion follows from the proof of Case 2.5.2 (B.1).

(D). There exists a $1\le i\le 4$, say $i=1$, such that  $G[x_{11},x_{12},x_{13},x_{14}]$ contains six edges.If $d_{G}(x_{i1})\leq k-2$, similar as Case 2.5.2 (B.1), we can prove the assertion.%let $G^{\prime}=G-\{u_{1},x_{i1}\}-\{x_{i2}x_{i3},x_{i2}x_{i4},x_{i3}x_{i4}\}$,$T^{\prime}=T-\{a_{0},a_{1}\}$.So $e(G^{\prime})\geq e(G)-(k-1)-(k-2)-3>\frac{1}{2}(k+4)(k-2)-(k-1)-(k-2)-3=\frac{1}{2}(k^{2}-2k-8)$.$\emph{avedeg}(G^{\prime})>(k^{2}-2k-8)/(k+2)=k-4$ and $\mid V(T^{\prime})\mid\leq k-4$.By the induction assumption,$T^{\prime}\subseteq G^{\prime}$.If $u_{1}$ hits $f^{\prime}(a_{2})$,then $f(a_{1})=u_{i}$.
%If $u_{i}$ misses $f^{\prime}(a_{2}),(a)f^{\prime}(a_{2})=x_{i2},f(a_{1})=x_{i1}$,if $ x_{i3}$ or $x_{i4} \notin f^{\prime}(V(T^{\prime}))$,then $f(a_{0})=x_{i3}$ or $x_{i4}$.else both of $x_{i3}$ and $x_{i4} \in f^{\prime}(V(T^{\prime})),u_{1}$ hits all the neighbours of $f^{\prime}(w)$,where $f^{\prime}(w)=x_{i3}$ or $x_{i4}$,
%so we let $f^{\prime}(w)=u_{1},f(a_{0})=x_{i3}$ or $x_{i4}$.(b)$f^{\prime}(a_{2})=x_{i3}$,the same as (a),(c)$f^{\prime}(a_{2})=x_{i4}$,the same as (a).All the situations are T-extensible.
 So we can assume $d_{G}(x_{i1})=d_{G}(x_{i2})=d_{G}(x_{i3})=d_{G}(x_{i4})=k-1$,we can also assume if $d_{G}(x)=k-1$,and $x$ misses $y$ then $d_{G}(y)=k-1$,furthermore we can assume $x$ hits all of the vertices whose degree is less tan $k-1$.
 let $G^{\prime}=G-\{u_{1},z\}$,z hits all of $\{x_{1},x_{2},x_{3},x_{4}\}$,$T^{\prime}=T-\{a_{0},a_{1}\}$.So $e(G^{\prime})\geq e(G)-(k-1)-(k-5)+1>\frac{1}{2}(k+4)(k-2)-(k-1)-(k-5)+1=\frac{1}{2}(k^{2}-2k+6)$.$\emph{avedeg}(G^{\prime})>(k^{2}-2k+6)/(k+2)>k-4$ and $\mid V(T^{\prime})\mid\leq k-2$.By the induction assumption,$T^{\prime}\subseteq G^{\prime}$.If $f^{\prime}(a_{2})$ hits $u_{1}$,then $f(a_{1})=u_{1},f(a_{0})=z$.If $f^{\prime}(a_{2})$ misses $u_{1}$,then $f(a_{0})=u_{1},f(a_{1})=z$.$f$ is T-extensible.

\end{document}